\documentclass{amsart}
\title{Model theory of valued fields}
\author[Raf Cluckers]{Raf Cluckers$^*$}
\thanks{$^*$Research Assistant of the Fund for Scientific Research --
 Flanders (Belgium)(F.W.O.)}
\address{Department of Mathematics\\
Katholieke Universiteit Leuven\\ Celestijnenlaan 200B\\ B-3001
Leuven, Belgium}
 \email{raf.cluckers@wis.kuleuven.ac.be}
 \urladdr{http://www.wis.kuleuven.ac.be/algebra/raf/}
\usepackage{amssymb}
\usepackage{amsmath}
\newtheorem{theorem}{Theorem}

\theoremstyle{definition}
\newtheorem{definition}{Definition}
\newtheorem{remark}{Remark}
\newtheorem{question}{Open Problem}
\newtheorem{example}{Example}

\DeclareMathOperator{\sq}{\square}
\begin{document}
\begin{abstract}
We give a proposal for future development of the model theory of
valued fields. We also summarize recent results on $p$-adic
numbers.
\end{abstract}
\maketitle
 \section{Valued fields} Let $K$ be a valued field
with a valuation map $v:K\to G\cup\{\infty\}$ to an ordered
group\footnote{Here, an ordered group is a totally ordered
non-trivial abelian group $G$  such that $x<y$ implies $x+z<y+z$
for all $x,y,z$ in $G$.} $G$; this is a map satisfying
\begin{itemize}
\item[(i)] $v(x)=\infty$ if and only if $x=0$;
\item[(ii)] $v(xy)=v(x)+v(y)$ for all $x,y\in K$;
\item[(iii)] $v(x+y)\geq \min\{v(x),v(y)\}$ for all $x,y\in K$.
\end{itemize}
We write $R$ for the  valuation ring $\{x\in K\mid v(x)\geq0\}$ of
$K$, $M$ for its unique maximal ideal $\{x\in K\mid v(x)>0\}$ and
we write $k$ for the residue field  $R/M$ and $\bar{}:R\to
k:x\mapsto\bar x$ for the natural projection. We call $K$ a
Henselian valued field if $R$ is a Henselian valuation ring.
\par
A valued field often carries an angular component map modulo $M$,
or angular component map for short; it is a group homomorphism
$\mathrm{ac}:K^\times\to k^\times$, extended by putting
$\mathrm{ac}(0)=0$, and satisfying $\mathrm{ac}(x)=\bar x$ for all
$x$ with $v(x)=0$ (see \cite{P}).

\section{Different languages}
 The model theory of valued fields can be studied at
different levels of complexity. One can use the most basic
language to study fields, the language $\mathcal{L}_{\rm ring}$ of
rings, and for some valued fields (like the $p$-adic numbers) the
relation $v(x)\leq v(y)$ is already definable in this language. In
general, this relation is not automatically definable (like in
algebraically closed valued fields). It is very natural to add a
divisibility predicate, or even more conveniently, a restricted
division function $D$ as follows
\[
D:K^2\to K^2:(x,y)\mapsto \left\{\begin{array}{ll} x/y  & \mbox{
if }
v(x)\geq v(y),\ y\not=0;\\
0 & \mbox{ else.}
\end{array}\right.
\]
Let us call $\mathcal{L}_D$ the language $\mathcal{L}_{\rm ring}$
together with $D$. Algebraically closed valued fields have
quantifier elimination in $\mathcal{L}_D$, see \cite{Rob_A}. It is
also convenient to add unary predicates $P_n$ to $\mathcal{L}_{\rm
ring}$, corresponding to the set of $n$-th powers in
$K^\times=K\setminus \{0\}$; one thus obtains the language of
Macintyre $\mathcal{L}_{\rm Mac}$. Macintyre \cite{Mac} proved
that $p$-adic fields have quantifier elimination\footnote{we
always allow parameters in formula's; thus by quantifier
elimination we mean that every formula is equivalent to a
quantifier free formula possibly containing parameters.} in
$\mathcal{L}_{\rm Mac}$. It is less known that there are many
other valued fields which have quantifier elimination in
$\mathcal{L}_{\rm Mac}$, like for example some fields of (formal)
Laurent series, $\mathbb{R}((t))$, $\mathbb{C}((t))$,
$\mathbb{Q}_p((t))$, and fields of iterated Laurent series like
$(\mathbb{R}((t_1)))((t_2))$, and so on. This follows from
quantifier elimination results of F.~Delon \cite{Delon}. Delon
adds more predicates (variants of the $P_n$ predicates) and
obtains quantifier elimination for a big class of valued fields;
this language $\mathcal{L}_{\rm Delon}$ reduces to
$\mathcal{L}_{\rm Mac}$ for the above mentioned fields. Let us
denote $\mathcal{L}_{{\rm Mac}, D}$ for the language
$\mathcal{L}_{\rm Mac}$ together with restricted division $D$.
\par
One can also add extra sorts to study valued fields, like the
value group, and the residue field. This way, one obtains
languages of Denef - Pas. More precisely, let $\mathcal{L}_k$ be
an arbitrary expansion of $\mathcal{L}_{\rm ring}$ and let
$\mathcal{L}_G$ be an arbitrary expansion of the language of
ordered groups with infinity $(+,-,0,\infty,\leq)$. A language of
Denef - Pas is a three-sorted language of the form
$(\mathcal{L}_k,\mathcal{L}_{\rm
ring},\mathcal{L}_G,v,\mathrm{ac})$ with as sorts:
\begin{itemize}
\item[(i)] a residue field-sort, with the language
$\mathcal{L}_k$,
 \item[(ii)] a valued field-sort, with the language $\mathcal{L}_{\rm
ring}$, and \item[(iii)] a value group-sort, with the language
$\mathcal{L}_G$.
\end{itemize}
The function symbol $v$ stands for the valuation map and
$\mathrm{ac}$ stands for an angular component map.
 \par
Pas obtains quantifier elimination for the valued-field-variables
for a very general class of fields, namely, for all Henselian
valued fields $K$ with an angular component map and with
$\mathrm{char} k=0$ (and also in other cases).
\par
The idea of these different languages is to be able to study more
and more complex kinds of valued fields, but of course also the
geometry of the definable sets gets more complex. In the following
sections we will define several notions of minimality with respect
to different languages.
\section{Basic open
problems} In this section, we let $K$ be a valued field which is
not real closed, nor algebraically closed. For the model theory of
real closed fields with a valuation, we refer to \cite{Dickmann}.
\begin{question}\label{question:Henselian}
%Let $K$ be a valued field.
Suppose that $K$ has quantifier elimination in $\mathcal{L}_{{\rm
Mac}, D}$, does it follow that $R$ is Henselian? If not, under
which extra conditions does this follow?
\end{question}
Similar question  can also be posed for languages of Denef - Pas
instead of $\mathcal{L}_{{\rm Mac}, D}$. One can put it even  more
challenging: does Henselianess of $R$ follow if we only know that
the definable subsets of $K$ are quantifier free definable in
$\mathcal{L}_{{\rm Mac}, D}$? Question \ref{question:Henselian} is
a generalization of an open problem mentioned in \cite{vdDMM},
where the analogue question is asked when $K$ is a $p$-valued
field.
\begin{question}[Same assumptions as in problem
\ref{question:Henselian}]\label{ques:indices} Does it follow that
the indices $[K^\times:P_n]$ are finite for $n>1$?
\end{question}
%Some other interesting question are:
\begin{question}[Same assumptions as
in problem \ref{question:Henselian}]\label{ques:not:alg:closed}
\begin{itemize}
\item
Does it follow that $R$ is already $\mathcal{L}_{\rm
Mac}$-definable, and hence, that $\mathcal{L}_{{\rm Mac}, D}$
is an expansion by definition of $\mathcal{L}_{\rm Mac}$?\\
\item Under which conditions does it follow that $G$ is a
$\mathbb{Z}$-group\footnote{By this we mean that $G$ is elementary
equivalent to $\mathbb{Z}$}, or that $G$ is elementary equivalent
to $\mathbb{Z}^n$ with lexicographical order\footnote{This is the
case for iterated Laurent series fields.}? \item Does it follow
that $K$ has characteristic zero?
\end{itemize}
 \end{question}
Also under the extra condition that the value group $G$ is a
$\mathbb{Z}$-group (or elementary equivalent to $\mathbb{Z}^n$
with lexicographical order) the above questions are worth looking
at. An example of a valued field with finite residue field and
value group $\mathbb{Z}^2$ with lexicographical order is
$\mathbb{Q}_p((t))$.
\section{A general notion of p-minimal fields}
In view of the above questions, %(possibly by requiring the condition that $G$ is elementary equivalent to $\mathbb{Z}$),
it is natural to give a definition of the following kind for $K$ a
valued field which is not algebraically closed, nor real closed.
\begin{definition}\label{def:p-minimal}
Let $\mathcal{L}$ be an expansion of $\mathcal{L}_{{\rm Mac}, D}$
and suppose that each $\mathcal{L}$-definable subset of $K$ is
already quantifier free definable
in $\mathcal{L}_{\rm Mac}$. %\footnote{Remark that this excludes algebraically closed valued fields and real closed fields.}
Then we call $(K,\mathcal{L})$ general p-minimal (where the p of
p-minimal comes from the powersets $P_n$). A theory of valued
fields in the language $\mathcal{L}$ is general p-minimal if each
of its models is.
\end{definition}
If $K$ is a $p$-adically closed field, this corresponds to the
notion of p-minimality by Haskell and Macpherson in
\cite{Haskell}. We can, of course, define analogously a notion of
Pas-minimality:
\begin{definition}\label{def:pas-minimal}
Let $\mathcal{L}_{\rm Pas}$ be a language of Denef - Pas and let
$\mathcal{L}$ be an expansion of $\mathcal{L}_{\rm Pas}$ and
suppose that each $\mathcal{L}$-definable subset of $K$ is already
quantifier free definable in $\mathcal{L}_{\rm Pas}$, then we call
$(K,\mathcal{L})$ Pas-minimal.
\end{definition}
On the $p$-adic numbers with either semialgebraic or subanalytic
structure, a cell decomposition theorem holds (see theorem
\ref{Thm:Cell} below, and also \cite{Denef} and \cite{Ccell}). In
a weaker formulation, Mourgues \cite{Mourgues} proves a cell
decomposition theorem for any p-minimal structure on
$\mathbb{Q}_p$ having definable Skolem functions. Mourgues'
decomposition is weaker in the sense that it only gives a
partition of definable sets into cells, and not a preparation of
definable functions (see theorem \ref{Thm:Cell} below for
preparations of definable functions). To obtain cell decomposition
\`a la Mourgues, definable Skolem functions are really needed. Not
all general p-minimal fields have definable Skolem functions, like
for example $\mathbb{C}((t))$; for p-minimal structures on
$\mathbb{Q}_p$ this is an open question. Languages of Denef - Pas
are very robust even if Skolem functions are not definable; the
residue field sort often replaces the need of definable Skolem
functions.
\par
It is natural to ask the following questions:
\begin{question}
\begin{itemize}
\item
Can one put on any general p-minimal field a p-minimal analytic
structure (as in section \ref{sec:p-adic} for $\mathbb{Q}_p$)
yielding
subanalytic sets\footnote{This is currently under development in \cite{CLipRob}.}?\\ %(The idea would be to add resticted analytic functions to the language as in \cite{DvdD}.)\\
\item Can one obtain cell decomposition (with preparation of
definable functions) for subanalytic
sets on a general p-minimal field with definable Skolem functions?\\
%\item
%Do p-minimal fields have definable Skolem functions?\\
\item
Do general p-minimal fields with definable Skolem functions have
 cell decomposition (with or without preparation)?
\end{itemize}
\end{question}
 Similar questions can be put for Pas - minimal fields and other
variants.
\section{Tame p-minimal fields}
\begin{definition}[Temporary definition in view of open problems
\ref{question:Henselian}, \ref{ques:indices}, and
\ref{ques:not:alg:closed}] Let $K$ be a valued field of
characteristic zero such that $G$ is elementary equivalent to
$\mathbb{Z}^n$ with lexicographical order. Let $\mathcal{L}$ be an
expansion of $\mathcal{L}_{{\rm Mac}, D}$. Suppose that
$(K,\mathcal{L})$ is a general p-minimal field, that $R$ is
Henselian and that the indices $[K^\times:P_n]$ are finite for
each $n>1$. Then we call $K$ a tame p-minimal field. We call a
theory in the language $\mathcal{L}$ tame p-minimal if each of its
models is.
\end{definition}
\begin{remark}
If the above problems \ref{question:Henselian}, \ref{ques:indices}
and \ref{ques:not:alg:closed} can be solved positively, any
general p-minimal field would be tame p-minimal.
\end{remark}
\begin{example}
An easy example of a tame p-minimal field is $\mathbb{R}((t))$,
the field of Laurent series over $\mathbb{R}$. If we fix a
parameter corresponding to $t$, then this field has definable
Skolem functions, because an ordering can be defined on
$\mathbb{R}((t))$ using cosets of the squares. Similar results
hold for $\mathbb{R}((t_1))((t_2))$ and so on. Other basic
examples are $\mathbb{Q}_p$, $\mathbb{C}((t))$, and
$\mathbb{Q}_p((t))$, and fields of iterated (formal) Laurent
series over $\mathbb{C}$, over $\mathbb{Q}_p$, and over finite
field extensions of $\mathbb{Q}_p$. This follows from the
quantifier elimination results of Delon.
\end{example}
In \cite{Cgroth} and \cite{CH}, criterions are given for the
existence of a $\mathcal{L}_{\rm ring}$-definable bijection $K\to
K^\times$ for a tame p-minimal field $K$; it is a sufficient
condition that the residue field $k$ is finite (like for
$\mathbb{Q}_p((t))$, $\mathbb{F}_p((t))((s))$ and so on).
\par
 I think it would be
very interesting to develop the model theory of tame p-minimal
fields and general p-minimal fields. For example, one can look for
a generalization of the notion of a $p$-adic cell (see section
\ref{sec:p-adic}) which can be used for all tame p-minimal fields
to obtain cell decomposition. The notion of $p$-adic cells itself
might be too strict: in general, the residue field is not finite
and therefore one might need cells of a more general kind.
Something in the style of the following subsets of $K$ might be
needed as cells:
 \begin{equation}%\label{Eq:cell}
\{x\in K\mid \wedge_{i,j,r} v(\alpha_{ir})\sq_{ir}
v(x-c_i)\lhd_{ij} v(x-c_j),\ x-c_j\in \lambda_j P_{n_j}\},
\end{equation}
where $\alpha_{ir}\in K^\times$, $c_j, \lambda_j\in K$, and
$\sq_{ir}$ and $\lhd_{ij}$ are either $<,>,=,$ or no condition,
and the conjunction is finite.
\begin{question}
Let $K$ be a tame p-minimal field. Does it follow that $K$ has
only finitely many algebraic field extensions of a given finite
degree?
\end{question}
\section{Fields of $p$-adic numbers}\label{sec:p-adic}
In this section we let $K$ be a finite field extension of
$\mathbb{Q}_p$ with valuation ring $R$ and $|\cdot|$ is the
$p$-adic norm. We put on $K$ the structure of subanalytic sets as
follows. For $X=(X_1,\ldots,X_m)$ let $R\langle X\rangle$ be the
ring of restricted power series over $R$ in the variables $X$; it
is the ring of power series $\sum a_iX^i$ in $R[[X]]$ such that
$|a_i|$ tends to $0$ as $|i|\to\infty$. (Here, we use the
multi-index notation where $i=(i_1,\ldots,i_m)$,
$|i|=i_1+\ldots+i_m$ and $X^i=X_1^{i_1}\ldots X_m^{i_m}$.) For
$x\in R^m$ and $f=\sum a_iX^i$ in $R\langle X\rangle$ the series
$\sum a_ix^i$ converges to a limit in $K$, thus, one can associate
to $f$ a \emph{restricted analytic function} given by
\[
f:K^m\to K:x\mapsto
 \left\{\begin{array}{ll} \sum_i a_i x^i & \mbox{ if
}x\in R^m,\\
0 & \mbox{ else.}
 \end{array}\right.
\]
We let $\mathcal{L}_{\rm sub}$ be the language $\mathcal{L}_{\rm
Mac}$ together with a function symbol $f$ for each $f\in R\langle
X\rangle$, with the interpretation as restricted analytic
function. A set $X\subset K^m$ is called (globally) subanalytic if
it is $\mathcal{L}_{\rm sub}$-definable. We call a function
$g:A\to B$ subanalytic if its graph is a subanalytic set. We refer
to \cite{DvdD}, \cite{Denef1}, \cite{vdDHM} and \cite{Ccell} for
the theory of subanalytic $p$-adic geometry and to \cite{Lip} for
the theory of rigid subanalytic sets.\\

$P$-adic cell decomposition makes use of basic sets called cells,
which we define inductively. For $m,l>0$ we write
$\pi_m:K^{m+l}\to K^m$ for the linear projection on the first $m$
variables, and, for $A\subset K^{m+l}$ and $x\in \pi_m(A)$, we
write $A_x$ for the fiber $\{t\in K^l\mid (x,t)\in A\}$.
 \\
Let $\mathcal{L}$ be an expansion of the language
$\mathcal{L}_{\rm Mac}$, then we define $\mathcal{L}$-cells.
 \begin{definition}\label{def:cell}
A cell $A\subset K$ is a (nonempty) set of the form
 \begin{equation}%\label{Eq:cell}
\{t\in K\mid |\alpha|\sq_1 |t-\gamma|\sq_2 |\beta|,\
  t-\gamma\in \lambda P_n\},
\end{equation}
with constants $n>0$, $\gamma,\lambda\in K$, $\alpha,\beta\in
K^\times$, and $\square_i$ either $<$ or no condition. If
$\lambda=0$ we call $A$ a (0)-cell and if $\lambda\not=0$, we call
$A$ a (1)-cell.
\par
 A $\mathcal{L}$-cell $A\subset K^{m+1}$ is a set  of the
form
 \begin{equation}\label{Eq:cell}
\{(x,t)\in K^{m+1}\mid x\in D, \  |\alpha(x)|\sq_1 |t-\gamma(
 x)|\sq_2 |\beta(x)|,\
  t-\gamma(x)\in \lambda P_n\},
\end{equation}
 with $(x,t)=(x_1,\ldots,
x_m,t)$, $n>0$, $\lambda\in K$, $D=\pi_m(A)$ a cell,
$\mathcal{L}$-definable functions $\alpha,\beta:K^m\to K^\times$
and $\gamma:K^m\to K$, and $\square_i$ either $<$ or no condition.
If $D$ is a $(i_1,\ldots,i_m)$-cell and $\lambda=0$ we call $A$ a
$(i_1,\ldots,i_m,0)$-cell, and if $\lambda\not=0$ we call $A$ a
$(i_1,\ldots,i_m,1)$-cell. Further, we call $\gamma$ the center of
the cell $A$ and we call $\lambda P_n$ the coset of $A$.
 \end{definition}
 \begin{remark}
 \begin{itemize}
 \item If $\mathcal{L}$ is $\mathcal{L}_{\rm Mac}$, we speak of semialgebraic cells; if
 $\mathcal{L}$ is $\mathcal{L}_{\rm sub}$, we speak of subanalytic cells.
 By induction we can call  $A$ an analytic cell if all
functions $\alpha,\beta,\gamma$ are analytic on $D$ and $D$ is an
analytic cell.
 \item An analytic $(i_1,\ldots,i_m)$-cell is a
$(\sum_j i_j)$-dimensional $p$-adic manifold. \item A cell
$A\subset K_1^{m+1}$ is either the graph of a subanalytic function
defined on $\pi_m(A)$ (if $\lambda=0$), or, for each $x\in
\pi_m(A)$, the fiber $A_x\subset K$ contains a nonempty open.
 \end{itemize}
\end{remark}
\par
Theorem \ref{Thm:Cell} below gives a cell decomposition for
subanalytic sets and for semialgebraic sets into analytic cells,
and at the same time it gives a preparation of definable
functions. In \cite{Denef1}, an overview is given of applications
of the semialgebraic cell decomposition.
 \begin{theorem}[$p$-adic cell decomposition]\label{Thm:Cell}
 Let $\mathcal{L}$ be either the language $\mathcal{L}_{\rm Mac}$ or $\mathcal{L}_{\rm sub}$.
Let $X\subset K^{m+1}$ be a $\mathcal{L}$-definable set, $m\geq0$,
and $f_j:X\to K$ $\mathcal{L}$-definable functions for
$j=1,\ldots,r$. Then there exists a finite partition of $X$ into
analytic $\mathcal{L}$-cells $A$ with center $\gamma:K^m\to K$ and
with coset $\lambda P_n$ such that for each $(x,t)\in A$
 \begin{equation}
 |f_j(x,t)|=
 |h_j(x)|\cdot|(t-\gamma(x))^{a_j}\lambda^{-a_j}|^\frac{1}{n},\qquad
 \mbox{ for each } j=1,\ldots,r,
 \end{equation}
with $(x,t)=(x_1,\ldots, x_m,t)$, integers $a_j$, and $h_j:K^m\to
K$ $\mathcal{L}$-definable functions, analytic on $\pi_m(A)$.  If
$\lambda=0$, we put $a_j=0$, and we use the convention $0^0=1$.
%Here we use the convention $0^0=1$ and if $\lambda=0$, $a_j=0$.
 \end{theorem}
The semialgebraic case is due to Denef \cite{Denef} and is stated
in this form in \cite{C}. The subanalytic case is recently proven
by the author in \cite{Ccell}. The fact that we can take the cells
and the functions $h_j$ to be analytic follows from the fact that
all semialgebraic and subanalytic functions are piecewise
analytic, see \cite{DvdD}. Namely, this fact can be used to
partition a given cell $A$ further into analytic cells with the
same center on which the restrictions of $h_j$ are analytic.

The proof of Thm.~ \ref{Thm:Cell} when
$\mathcal{L}=\mathcal{L}_{\rm sub}$ uses compactness arguments and
recent results of \cite{vdDHM}. Once one knows this cell
decomposition, it follows that $(K,\mathcal{L}_{\rm sub})$ is
p-minimal (although the subanalytic cell decomposition actually
relies on this fact, proven in \cite{vdDHM}). Recently, it has
been proven by D.~Haskell and the author (unpublished notes) that
if one adds a nontrivial entire analytic function to
$\mathcal{L}_{\rm Mac}$, the obtained structure on $\mathbb{Q}_p$
is not p-minimal anymore.
\par
An important use of $p$-adic cell decomposition is to calculate
$p$-adic integrals, see e.g.~\cite{Denef3}. We explain this into
some detail.
\begin{definition}\label{basic algebra's}
Let $\mathcal{L}$ be either the language $\mathcal{L}_{\rm Mac}$
or $\mathcal{L}_{\rm sub}$. The algebra
$\mathcal{C}_{\mathcal{L}}(K^m)$ of ($\mathcal{L}$-definable)
constructible functions on $K^m$ is the $\mathbb{Q}$-algebra
generated by the functions $K^m\to\mathbb{Q}$ of the form
$x\mapsto v(h(x))$ and $x\mapsto |h'(x)|$ where $h:K^m\to
K^\times$ and $h':K^m\to K$ are $\mathcal{L}$-definable.
\par
To any function $f\in \mathcal{C}_{\mathcal{L}}(K^{m+n})$ ,
$m,n\geq 0$, we associate a function $I_m(f):K^m\to \mathbb{Q}$ by
putting
 \begin{equation}\label{I_l}
I_m(f)(\lambda)= \int\limits_{K^n}f(\lambda,x)|dx|
 \end{equation}
if the function $x\mapsto f(\lambda,x)$ is absolutely integrable
for all $\lambda\in K^m$, and by putting $I_m(f)(\lambda)=0$
otherwise. Of course, $|dx|$ stands here for a Haar measure.
\end{definition}
Constructible functions often appear naturally in number theory,
for example, local singular series, introduced by A.~Weil, are
constructible (see \cite{Denef1} and \cite{Ccell}).
 \begin{theorem}[Basic Theorem on $p$-adic Analytic Integrals]\label{thm:basic}
For any function $f\in \mathcal{C}_{\mathcal{L}}(K^{m+n})$, the
function $I_m(f)$ is in $\mathcal{C}_{\mathcal{L}}(K^m)$.
 \end{theorem}
The semialgebraic case of Thm.~\ref{thm:basic} is proven in
\cite{Denef1}, the subanalytic case in \cite{Ccell}. The following
questions are natural to ask.
\begin{question}
\begin{itemize}
\item Does any p-minimal structure on $\mathbb{Q}_p$ have definable Skolem
functions?
\item Does any p-minimal structure on $\mathbb{Q}_p$ allow a cell
decomposition with preparation of definable functions?
\item Apart from the semialgebraic and subanalytic structures on
$\mathbb{Q}_p$, can one find other p-minimal structures on
$\mathbb{Q}_p$?
\end{itemize}
\end{question}
\bibliography{cluckers}

\begin{thebibliography}{10}

\bibitem{C}
R.~Cluckers.
\newblock Classification of semialgebraic $p$-adic sets up to semialgebraic
  bijection.
\newblock {\em Journal f{\"u}r die reine und angewandte Mathematik},
  540:105--114, 2001.

\bibitem{Ccell}
R.~Cluckers.
\newblock Analytic $p$-adic cell decomposition and integrals.
\newblock {\em Trans. Am. Math. Soc.}, To Appear.
\newblock Preprint available at arXiv:math.NT/0206161.

\bibitem{Cgroth}
R.~Cluckers.
\newblock Grothendieck rings of laurent series fields.
\newblock {\em Journal of Algebra}, to appear.
\newblock available at arXiv:math.LO/0210350.

\bibitem{CH}
R.~Cluckers and D.~Haskell.
\newblock {G}rothendieck rings of $\mathbb{Z}$-valued fields.
\newblock {\em Bulletin of Symbolic Logic}, 7(2):262--269, 2001.

\bibitem{CLipRob}
R.~Cluckers, L.~Lipshitz, and Z.~Robinson.
\newblock Analytic cell decomposition.
\newblock Preprint.

\bibitem{Delon}
F.~Delon.
\newblock {\em Quelques propri\'et\'es des corps valu\'es}.
\newblock th\`ese d'\'etat, Universit\'e Paris VII, 1981.

\bibitem{Denef}
J.~Denef.
\newblock The rationality of the {P}oincar\'e series associated to the $p$-adic
  points on a variety.
\newblock {\em Inventiones Mathematicae}, 77:1--23, 1984.

\bibitem{Denef3}
J.~Denef.
\newblock On the evaluation of certain $p$-adic integrals.
\newblock In {\em Th\'eorie des nombres, S\'emin. Delange-Pisot-Poitou
  1983--84}, volume~59, pages 25--47, 1985.

\bibitem{Denef1}
J.~Denef.
\newblock {\em Arithmetic and geometric applications of quantifier elimination
  for valued fields}, volume~39 of {\em MSRI Publications}, pages 173--198.
\newblock Cambridge University Press, 2000.

\bibitem{DvdD}
J.~Denef and L.~van~den Dries.
\newblock $p$-adic and real subanalytic sets.
\newblock {\em Annals of Mathematics}, 128(1):79--138, 1988.

\bibitem{Dickmann}
M.~A. Dickmann.
\newblock Elimination of quantifiers for ordered valuation rings.
\newblock {\em Journal of Symboloc Logic}, 52:116--128, 1987.

\bibitem{Haskell}
D.~Haskell and D.~Macpherson.
\newblock A version of o-minimality for the $p$-adics.
\newblock {\em Journal of Symbolic Logic}, 62(4):1075--1092, 1997.

\bibitem{vdDHM}
D.~Haskell, D.~Macpherson, and L.~van~den Dries.
\newblock One-dimensional $p$-adic subanalytic sets.
\newblock {\em Journal of the London Mathematical Society}, 59(1):1--20, 1999.

\bibitem{Lip}
L.~Lipshitz and Z.~Robinson.
\newblock {\em Rings of separated power series and quasi-affinoid geometry},
  volume 264 of {\em Paris: Soci\'et\'e Math\'ematique de France}.
\newblock Ast\'erisque, 2000.

\bibitem{Mac}
A.~Macintyre.
\newblock On definable subsets of $p$-adic fields.
\newblock {\em Journal of Symbolic Logic}, 41:605--610, 1976.

\bibitem{vdDMM}
A.~Macintyre, K.~McKenna, and L.~van~den Dries.
\newblock Elimination of quantifiers in algebraic structures.
\newblock {\em Advances in mathematics}, 47:74--87, 1983.

\bibitem{Mourgues}
M.-H. Mourgues.
\newblock Corps p-minimaux avec fonctions de {S}kolem d{\'e}finissables.
\newblock pages 1--8, 1999.

\bibitem{P}
J.~Pas.
\newblock On the angular component map modulo $p$.
\newblock {\em Journal of Symbolic Logic}, 55:1125--1129, 1990.

\bibitem{Rob_A}
Abraham Robinson.
\newblock {\em Complete theories}.
\newblock Studies in Logic and the Foundations of Mathematics. North- Holland
  Publishing, Amsterdam, 1956.
\newblock 129 p.

\end{thebibliography}
\bibliographystyle{plain}
\end{document}